\documentclass[conference]{IEEEtran}

\usepackage{microtype}
\usepackage{graphicx}
\usepackage{subfigure}
\usepackage{amsmath,amssymb}
\usepackage{booktabs} 

\usepackage{hyperref}

\newtheorem{theorem}{Theorem}
\newtheorem{lemma}{Lemma}

\begin{document}

\title{Chebyshev Inertial Iteration for Accelerating Fixed-Point Iterations}

\author{%
  \IEEEauthorblockN{
  		Tadashi Wadayama\IEEEauthorrefmark{1}
  		and Satoshi Takabe\IEEEauthorrefmark{1}\IEEEauthorrefmark{2}}
  \IEEEauthorblockA{\IEEEauthorrefmark{1}%
		Nagoya Institute of Technology,
		Gokiso, Nagoya, Aichi, 466-8555, Japan,\\
 		\{wadayama, s\_takabe\}@nitech.ac.jp} 
 \IEEEauthorblockA{\IEEEauthorrefmark{2}%
		RIKEN Center for Advanced Intelligence Project, Chuo-ku, Tokyo, 103-0027, Japan
 		} \\ 
}

\maketitle

\begin{abstract}
A novel method which is called the Chebyshev inertial iteration for accelerating the convergence speed of 
fixed-point iterations is presented. 
The Chebyshev inertial iteration can be regarded as a valiant of
the successive over relaxation or Krasnosel'ski\v{\i}-Mann iteration
utilizing the inverse of roots of a Chebyshev polynomial as 
iteration dependent inertial factors.
One of the most notable features of the proposed method is that it can be applied to
nonlinear fixed-point iterations in addition to linear fixed-point iterations.
Linearization around the fixed point is the key 
for the analysis on the local convergence rate of the proposed method.
The proposed method appears effective in particular for accelerating the proximal gradient methods 
such as ISTA. It is also proved that the proposed method can successfully accelerate almost any fixed-point iterations if all the eigenvalues of 
the Jacobian at the fixed point are real. 
\end{abstract}

\section{Introduction}
\label{submission}

A {\em fixed-point iteration} 
\begin{equation}
x^{(k+1)} = f(x^{(k)}), \ k = 0, 1, 2, \ldots 
\end{equation}
 is universally used for 
solving scientific and engineering problems
such as inverse problems \cite{Bauschke11a}. 
 An eminent example is 
a gradient descent method for minimizing a convex or non-convex objective function \cite{Boyed04}.
Variants of gradient descent methods such as stochastic gradient descent methods  
is becoming indispensable for solving large scale machine 
learning problems. The fixed-point iteration is also widely employed for
solving large scale linear equations.
For example, Gauss-Seidel methods, Jacobi methods, and 
conjugate gradient methods  are
well known fixed-point iterations to solve linear equations \cite{Saad03}.

Another example of fixed-point iterations is a proximal gradient descent method for solving 
a certain class of convex problems. A notable instance is Iterative Shrinkage-Thresholding Algorithm (ISTA) \cite{Daubechies04}
for sparse signal recovery problems. The projected gradient method is also included  in the class of
the proximal gradient method.

In this paper, we present a novel method for accelerating the convergence speed of
linear and nonlinear fixed-point iterations.  The proposed method is called 
the {\em Chebyshev inertial iteration} because it is
heavily based on the property of the Chebyshev polynomials \cite{Mason03}.
The proposed acceleration method is closely related to successive over relaxation (SOR) as well.
The SOR is well known as a method for accelerating Gauss-Seidel and Jacobi method to solve
a linear equation $P x = q$ where $P \in \mathbb{R}^{n \times n}$ \cite{Young71} .
For example, the Jacobi method is  based on the following linear fixed-point iteration: 
\begin{equation}
	x^{(k+1)} = f(x^{(k)}) := D^{-1} \left(q - (P - D) x^{(k)}  \right),
\end{equation}
where $D$ is the diagonal matrix whose diagonal elements are identical to the diagonal elements of $P$.
The SOR for the Jacobi iteration is simply given by
\begin{equation}
x^{(k+1)} = x^{(k)} + \omega_k\left( f( x^{(k)} ) - x^{(k)}  \right), 
\end{equation}
where $\omega_k \in (0, 2)$. 
A fixed SOR factor $\omega_k := \omega$ is often employed in practice. 
When the SOR factors $\{\omega_k \}$ are appropriately chosen, the SOR-Jacobi method 
achieves much faster convergence compared with the plain Jacobi method.

\begin{figure}[tbp]
\vskip 0.2in
\begin{center}
\centerline{\includegraphics[width=\columnwidth]{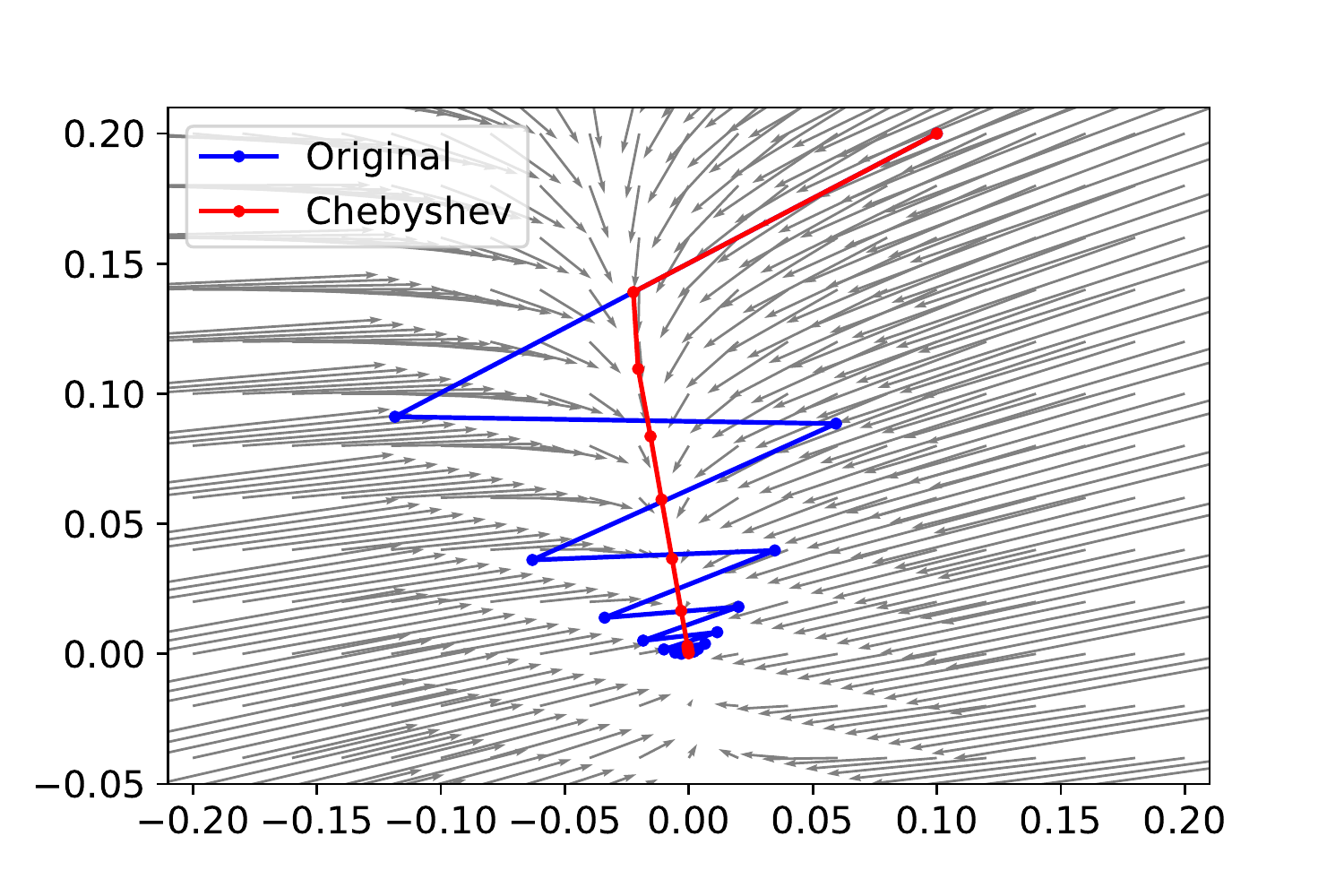}}
\caption{Trajectories for $x^{(k+1)} = \tanh(A x^{(k)})$: The zigzag line (blue) is the trajectory of the original fixed-point iteration. The fixed point is $(0, 0)^T$. Another line (red) is the trajectory of the Chebyshev inertial iteration with $T=8$. 
}
\label{cover}
\end{center}
\vskip -0.2in
\end{figure}

The Chebyshev inertial iteration proposed in this paper can be regarded as a valiant of
the SOR method or Krasnosel'ski\v{\i}-Mann iteration \cite{Bauschke11b} utilizing the inverse of roots of a Chebyshev polynomial as {\em iteration dependent inertial factors}.
This choice reduces the spectral radius of matrices related 
to the convergence of the fixed-point iteration.
One of the most notable features of the proposed method is that it can be applied to
{\em nonlinear} fixed-point iterations in addition to linear fixed-point iterations.
Linearization around the fixed point is the key for the analysis on
 the local convergence rate of the 
Chebyshev inertial iteration.

Figure \ref{cover} indicates an example of the Chebyshev inertial iteration for 
a two-dimensional nonlinear fixed-point iteration
 $x^{(k+1)} = \tanh(A x^{(k)})$.
The initial point is $(0.1, 0.2)^T$. The matrix $A$ is defined by
$
	A := 
	\left(
	\begin{array}{cc}
-0.6929 & -0.2487 \\
-0.2870 &  0.6005 \\
	\end{array}
	\right).
$
 The length of background arrows is proportional to $\tanh(A x) - x$.
We can observe that the trajectory of 
the Chebyshev inertial iteration (red)  straightly reaches the fixed point 
with fewer iterations.

In this paper, it will be proved that the proposed method can successfully accelerate 
almost any fixed-point iterations if all the eigenvalues of 
the Jacobian at the fixed point are real.

\section{Cyebyshev Inertial Iterations}

\subsection{Inertial iteration}

Let $f: \mathbb{R}^n \rightarrow \mathbb{R}^n$ be a differentiable function.
Throughout the paper, we will assume that $f$ has a fixed point $x^*$ satisfying 
$x^* = f(x^*)$ and that $f$ is a {\em locally contracting mapping} around $x^*$.
Namely, there exists $\epsilon (\epsilon > 0)$ such that
$
	\|f(x) - f(y) \| < \|x - y\|
$
holds for any $x, y (x \ne y) \in B(x^*, \epsilon)$ where
\begin{equation}
B(x^*, \epsilon) := \{x \in \mathbb{R}^n: \|x - x^*\| \le \epsilon \}.	
\end{equation}
The norm $\|\cdot \|$ represents Euclidean norm.
From these assumptions, it is clear that the fixed-point iteration
\begin{equation} \label{iteration}
x^{(k+1)} = f(x^{(k)}),\   k = 0, 1, 2, \ldots
\end{equation}
eventually converges to the fixed point $x^*$  if the initial point 
$x_0$ is included in $B(x^*, \epsilon)$.

The SOR is a well known technique for accelerating the 
convergence speed for linear fixed point iteration such as Gauss-Seidel method \cite{Young71}.
The SOR (or Krasnosel'ski\v{\i}-Mann) iteration can be generalized into the {\em inertial iteration} for fixed-point iterations, which is given by
\begin{eqnarray} \nonumber
	x^{(k+1)} &=& x^{(k)} + \omega_k \left[f(x^{(k)}) 
	- x^{(k)} \right] \\ \label{SORiteration}
	&=& (1 - \omega_k)x^{(k)} + \omega_k f(x^{(k)}),
\end{eqnarray}
where $\omega_k (k = 0, 1, \ldots)$ is a real coefficient, which is
called an {\em inertial factor}.  When $f$ is linear and $\omega_k  = \omega$ is a constant, 
the optimal choice of $\omega$ in terms of the convergence rate is known \cite{Saad03} but
the optimization of iteration dependent inertial factors ${\omega_k}$  for  nonlinear update functions
has not been studied as far as the authors know.

We here define $S^{(k)}: \mathbb{R}^n \rightarrow \mathbb{R}^n$ by
\begin{equation}
S^{(k)}(x) := (1 - \omega_k) x + \omega_k f(x).
\end{equation}
It should be remarked that 
\begin{equation}
	S^{(k)}(x^*) = (1 - \omega_k) x^{*} + \omega_k f(x^{*}) = x^*
\end{equation}
holds for any $k$.
Since  $x^*$ is also the fixed point of $S^{(k)}$, 
the inertial iteration
does not change the fixed point of the original fixed-point iteration.

\subsection{Spectral condition}

Since $f$ is a differentiable function, 
it would be natural to consider the {\em linear approximation} around the fixed point $x^*$. 
The Jacobian matrix of $S^{(k)}(x)$ around $x^*$
 is given by $I - \omega_k (I- J^*)$ where $I$ represents the identity matrix.
The matrix $J^*$ is the Jacobian matrix of $f$ at  $x^*$:
\begin{equation}
	J^* := 
	\left(
	\begin{array}{ccc}
		\frac{\partial f_1}{\partial x_1}(x_1^*) &\frac{\partial f_1}{\partial x_2}(x_2^*) & \hdots \\	
		 \vdots & \ddots & \vdots \\	
		\frac{\partial f_n}{\partial x_1}(x_1^*) &\hdots & \frac{\partial f_n}{\partial x_n}(x_n^*) \\	
	\end{array}
	\right),
\end{equation}
where $f := (f_1, f_2, \ldots, f_n)^T$ and $x^* := (x_1^*,\ldots, x_n^*)^T$.
The convergence rate of the inertial iteration is dominated by the 
eigenvalues of $I - J^*$.
In the following discussion, we will use the abbreviation
$
	B := I - J^*.
$
The next lemma provides useful information on the range of the eigenvalues of $B$
whose proof is given in Appendix.
\begin{lemma}
	Assume that the function $f$ is a locally contracting mapping around $x^*$ 
	and that all the eigenvalues of $J^*$ are real.
	All eigenvalues of $B$ lie in the open interval  $(0,  2)$.
\end{lemma}

Let us get back to the discussion on linearization around the fixed point.
The function $S^{(k)}(x)$ in the inertial iteration 
can be approximated by
\begin{equation}
	S^{(k)}(x) \simeq S^{(k)}(x^*) +  \left(I - \omega_{k} B \right) (x - x^*) 
\end{equation}
due to Taylor series expansion around the fixed point\footnote{In this context,  we will omit the residual terms such as $o(\|x - x^*\|)$ appearing in the
Taylor series to simplify the argument. More careful treatment on the residual terms 
can be found in Appendix.}.
By letting $x = x^{(k)}$, we have
\begin{equation} \label{aprox}
	x^{(k+1)} - x^* \simeq  \left(I - \omega_{k} B \right) (x^{(k)} - x^*).
\end{equation}
Applying a norm inequality to  (\ref{aprox}), we immediately obtain
\begin{equation} 
	\|x^{(T)} - x^* \|  \lesssim \rho\left(\prod_{k = 0}^{T-1} \left(I - \omega_{k} B \right) \right)  \|x^{(0)} - x^*\|,
\end{equation}
where $\rho(\cdot)$ represents the spectral radius.

In the following argument, we assume 
periodical inertial factors $\{\omega_k \}$ satisfying 
$\omega_{\ell T + j} = \omega_j$ ($\ell = 0,1,2,\ldots$, $j = 0, 1, 2, \ldots, T-1$)
where $T$ is a positive integer.
If the spectral radius satisfies 
\begin{equation} \label{spectral}
\rho_T := \rho\left(\prod_{k = 0}^{T-1} \left(I - \omega_{k} B \right) \right) < 1,
\end{equation}
then we can expect local linear convergence around the fixed point:
\begin{equation} \label{lT}
	\|x^{(\ell T)} - x^* \|  \lesssim \rho_T^\ell  \|x^{(0)} - x^*\|.
\end{equation}
for positive integer  $\ell$.
The important problem is to find an appropriate set of inertial factors 
satisfying the {\em spectral condition} (\ref{spectral}) to ensure the convergence.

\subsection{Chebyshev inertial iteration}

As described in the previous subsection, the choice of inertial factors $\{\omega_k\}$ is crucial to obtain 
local convergence of the inertial iteration.
However, the direct minimization of the spectral radius $\rho_T$ seems computationally intractable 
because $\rho_T$ is nonconvex with respect to $\{\omega_k\}$.
In this subsection, we will show that the {\em Chebyshev inertial factors} defined by
\begin{equation} \label{SORfactors}
\omega_k
    :=\left[ \lambda_+ + \lambda_- \cos \left(\frac{2k+1}{2T} \pi \right) \right]^{-1}, \  k = 0, \ldots, T-1,
\end{equation}
where
\begin{eqnarray}
\lambda_+ &:=& \frac{\lambda_{max}(B) + \lambda_{min}(B)}{2}, \\
\lambda_- &:=& \frac{\lambda_{max}(B) - \lambda_{min}(B)}{2},
\end{eqnarray}
satisfies the spectral condition (\ref{spectral}) 
and can significantly improve the local convergence rate.
The $\lambda_{min}(\cdot)$ and $\lambda_{max}(\cdot)$ represents
the minimum and maximum eigenvalues, respectively.
The inertial iteration using the Chebyshev inertial factors is referred to as
the {\em Chebyshev inertial iteration} hereafter.

The Chebyshev polynomial $C_k(x)$ is recursively defined by
\begin{equation}
	C_{k+1}(x)	 = 2 x C_k(x) - C_{k-1}(x)	
\end{equation}
with the initial conditions $C_0(x) = 1$ and $C_1(x) = x$ \cite{Mason03}. 
One of the significant properties of the Chebyshev polynomial is the minimality 
in $\infty$-norm of polynomials defined on the closed interval $[-1, 1]$.
Namely, the monic Chebyshev polynomial $2^{1-T} C_T(x)$,
i.e., the coefficient of leading term is normalized to one,  gives the minimum $\infty$-norm
\begin{equation}
	\max_{-1 \le x \le 1}|2^{1-T} C_T(x)| = 2^{1-T}
\end{equation}
among any monic polynomials of degree $T$.
In the following, we consider the {\em affine transformed monic Chebyshev (ATMC) polynomial}
on the closed interval $[a, b] (a \neq b)$:
\begin{equation}
	\hat C_T(x; a, b):= 2^{1-T} \left(\frac{b - a}{2}\right)^T 
	C_T\left(\frac{2 x - b - a}{b- a} \right).
\end{equation}
The roots $\{z_0, z_1,\ldots, z_{T-1}\}$ of the ATMC polynomial $\hat C_T(x; a, b)$ are given by
\begin{equation}
	z_k := \frac{b + a}{2} + \frac{b-a}{2}\cos \left(\frac{2k + 1}{2 T} \pi   \right),\ k = 0,  \ldots, T-1,
\end{equation}
which are called the {\em Chebyshev roots}.

The matrix polynomial $\prod_{k = 0}^{T-1} \left(I - \omega_{k} B \right)$ in the spectral condition (\ref{spectral})
corresponds to a polynomial defined on $\mathbb{R}$:
\begin{equation}
	\beta_T(\lambda) := \prod_{k = 0}^{T-1} \left(1 - \omega_{k} \lambda \right)
	=  \left(\prod_{k = 0}^{T-1} \omega_k \right)\prod_{k = 0}^{T-1} \left(\frac{1}{\omega_{k}} -  \lambda \right).
\end{equation}
The following analysis is based on the fact that
any eigenvalue of $\beta_T(B) = \prod_{k = 0}^{T-1} \left(I - \omega_{k} B \right)$ 
can be expressed as $\beta_T(\lambda')$
where $\lambda'$ is an eigenvalue of $B$ \cite{Courant53}.

\begin{figure}[t]
\vskip 0.2in
\begin{center}
\centerline{\includegraphics[width=\columnwidth]{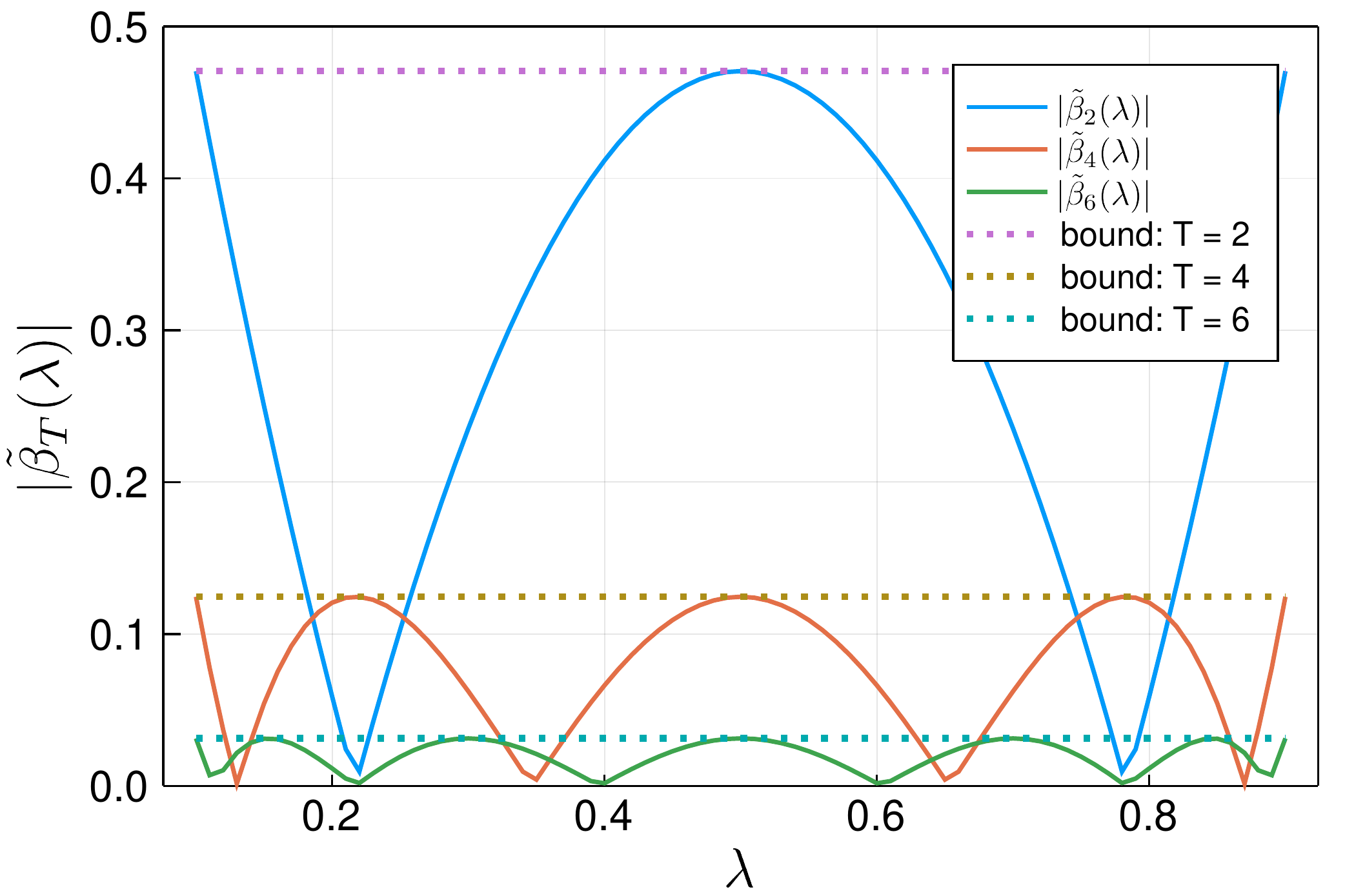}}
\caption{The absolute value of $\tilde \beta_T(\lambda) (T = 2, 4, 6)$:  The minimum and maximum eigenvalues of $B$ is assumed to be $\lambda_{min}(B) = 0.1$
	and $\lambda_{max}(B) = 0.9$. The dashed lines represent the values of  the upper bound
	$\text{\rm sech}\left(T \cosh^{-1}\left(\frac{0.9 + 0.1}{0.9 - 0.1} \right)\right)$  for $T = 2, 4, 6$.}
\label{monic}
\end{center}
\vskip -0.2in
\end{figure}

Our strategy to choose the inertial factor is {\em to let  $\omega_k$ be the inverse of 
the Chebyshev root}, i.e, 
$
\omega_k = \tilde \omega_k := {1}/{z_k}
$
in order to
decrease the absolute value of $\beta_T(\lambda)$ under the assumption
$a := \lambda_{min}(B)$ and $b := \lambda_{max}(B)$.
In other words, we embed the ATMC polynomial $\hat C_T(\lambda; a, b)$ into $\beta_T(\lambda)$
with expectation such that $\hat C_T(\lambda; a, b)$ can provide
small $|\beta_T(\lambda)|$  due to the minimality of the Chebyshev polynomials.
Let 
\begin{equation}
\tilde \beta_T(\lambda) := \prod_{k = 0}^{T-1} \left(1 - \tilde \omega_{k} \lambda \right).
\end{equation}
Since $\prod_{k = 0}^{T-1} \left(z_k -  \lambda \right) = \hat C_T(\lambda; a, b)$,
we have
\begin{eqnarray}
\tilde \beta_T(\lambda) &=& \left(\prod_{k = 0}^{T-1} \frac{1}{z_k} \right)	
\hat C_T(\lambda; a, b) \\ 
&=& \left[\hat C_T(0; a, b) \right]^{-1}
\hat C_T(\lambda; a, b),
\end{eqnarray}
where $z_1 z_2 \cdots z_{T-1} = \hat C_T(0; a, b)$.

The following lemma can be immediately derived based on a property on the Chebyshev polynomials.
\begin{lemma}
Let $\omega_k = \tilde \omega_k$ for $k = 0, 1, \ldots, T-1$.
The absolute value of $\tilde \beta_T(\lambda)$ can be upper bounded by
\begin{equation}
|\tilde \beta_T(\lambda)| \le  \text{\rm sech}\left(T \cosh^{-1}\left(\frac{b+a}{b-a}\right) \right)
\end{equation}
for $\lambda \in [a, b]$.
\end{lemma}
The proof of the lemma is given in Appendix.

Figure \ref{monic} presents the plot of the absolute values of $\tilde \beta_T(\lambda)$ for $T = 2, 4, 6$
under the assumption  $a = \lambda_{min}(B) = 0.1$ and $b = \lambda_{max}(B) = 0.9$.

It is clear that $\cosh^{-1}\left((b+a)/(b-a)\right) > 0$ holds.
We thus have
\begin{equation}
	\rho_T \le  \max_{a \le \lambda \le b}|\tilde \beta_T(\lambda)| = \text{sech}\left(T \cosh^{-1}\left(\frac{b+a}{b-a}\right) \right) < 1
\end{equation}
because $0 < \text{sech}(\alpha x) < 1$ holds for $x > 0$ and $\alpha > 0$.
This means that the spectral condition is satisfied  when 
the inertial factors are set to the inverse of the Chebyshev roots.

These arguments can be summarized as follows.
\begin{theorem} \label{convergence_rate}
Assume that all the eigenvalues of $J^*$ are real.
Let $x^{(0)}, x^{(1)},\ldots$ be the point sequence obtained by the Chebyshev inertial iteration
(\ref{SORiteration}). Let $\delta = C \|x_0 - x^*\|$ where 
$C$ is a positive constant.
For $\ell = 1,2,3, \ldots$, we have
\begin{equation}
	\lim_{\delta \rightarrow 0}\frac{\|x^{(\ell T)} - x^*\|}{\|x^{(0)} - x^*\|}	
	\le  \left[\text{\rm sech}\left(T \cosh^{-1}\left(\frac{b+a}{b-a}\right) \right) \right]^\ell,
\end{equation}
where $a := \lambda_{min}(B)$ and $b := \lambda_{max}(B)$.
\end{theorem}

The proof of the theorem is given in Appendix.
The claim of the theorem provides an upper bound on the local convergence rate 
of the Chebyshev inertial iteration.

\subsection{Acceleration of convergence}

In the previous subsection, we saw the spectral radius $\rho_T$ is 
strictly smaller than one when the Chebyshev inertial factors are applied. 
We still need to confirm whether the Chebyshev inertial iteration
certainly accelerates the convergence compared with the convergence rate of
the original fixed-point iteration.

\begin{figure}[t]
\vskip 0.2in
\begin{center}
\centerline{\includegraphics[width=\columnwidth]{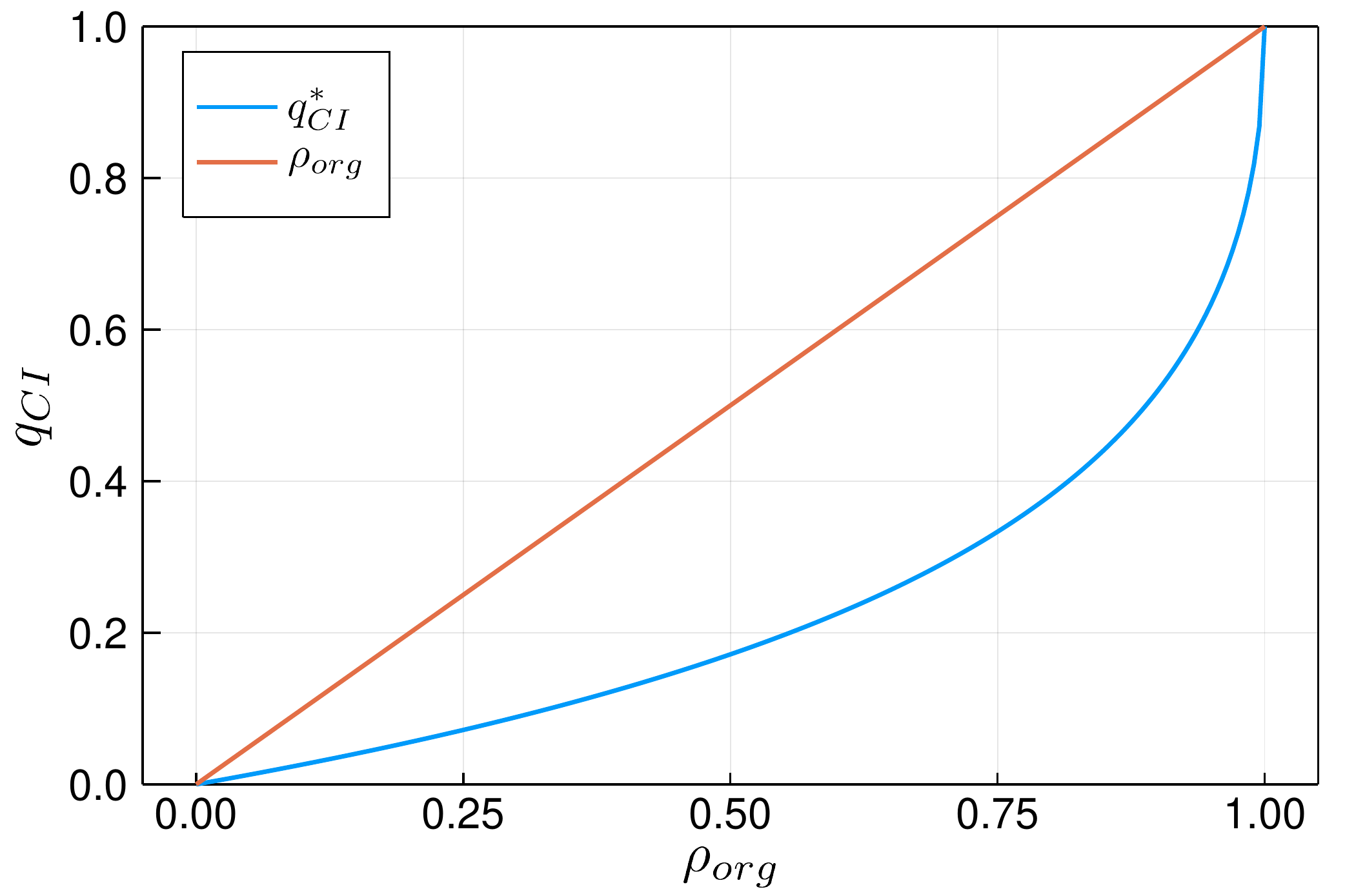}}
\caption{The convergence rate region of $q_{CI}(T)$: The straight line represents $\rho_{org}$ 
and the curve corresponds to $q_{CI}^*$.
For sufficiently large $T$, $q_{CI}(T)$ falls in the area surrounded by 
the line and the curve.}
\label{xi_function}
\end{center}
\vskip -0.2in
\end{figure}

Around the fixed point $x^*$, the original fixed-point iteration (\ref{iteration})
achieves 
$\|x^{(k)} - x^*\| \le q_{org}^k \|x^{(0)} - x^*\|$ where $q_{org} := \rho(J^*) < 1$.
In order to make fair comparison, we here define
\begin{equation}
q_{CI}(T):=\left[\text{\rm sech}\left(T \cosh^{-1}
\left(\frac{b+a}{b-a}\right) \right) \right]^{1/T}
\end{equation}
for the Chebyshev inertial iteration.
From Theorem \ref{convergence_rate}, 
\begin{equation}
\|x^{(k)} - x^*\| < q_{CI}^k \|x^{(0)} - x^*\|	
\end{equation}
holds if $k$ is a multiple of $T$ and $x^{(0)}$ is sufficiently close to $x^*$.
We can show that $q_{CI}(T)$ is a bounded monotonically decreasing function.
The limit value of $q_{CI}(T)$ is given by 
\begin{equation}
q_{CI}^* := \lim_{T \rightarrow \infty} q_{CI}(T) = 
\exp\left(-\cosh^{-1}\left(\frac{b+a}{b-a}\right)  \right),
\end{equation}
by using 
$
\lim_{x \rightarrow \infty}\left[\text{sech}(\alpha x) \right]^{1/x} = e^{-\alpha}.
$
In the following discussion, we will assume 
that $\lambda_{min}(J^*) = 0$ and $\lambda_{max}(J^*) = \rho_{org} (0 < \rho_{org} < 1)$.
In this case we have $a := 1- \rho_{org}$ and $b := 1$.
We can show that, for sufficiently large $T$, we have 
$
q_{CI}^* < q_{CI}(T) < \rho_{org}.
$
Figure \ref{xi_function} presents the region where $q_{CI}(T)$ can exist.
This implies that the Chebyshev inertial iteration can certainly achieve acceleration of the 
local convergence speed if $T$ is sufficiently large and $x_0$ is sufficiently close to $x^*$.

\section{Proximal Iteration}

In this section, we will discuss an important special case where 
$f$ has the form 
$
f(x) := g(A x + b)	
$
that is a composition of the affine transformation $A x + b$ and 
a component-wise nonlinear function $g$ where 
$A \in \mathbb{R}^{n \times n}, b \in \mathbb{R}^n$.
Note that a function $g:\mathbb{R} \rightarrow \mathbb{R}$ is applied to 
$x := (x_1,\ldots, x_n)^T \in \mathbb{R}^n$ in component-wise manner as
$
	g(x) := (g(x_1), \ldots, g(x_n))^T.
$

The fixed-point iteration defined by
\begin{equation} \label{affine}
x^{(k+1)} = g(A x^{(k)} + b), \ k = 0, 1, \ldots 
\end{equation}
is referred to as a {\em proximal iteration} 
because this type of iteration often appears in proximal gradient descent methods \cite{Parikh14}.
The proximal iteration is important because 
many iterative optimization methods such as projected gradient methods
are included in the class.

\subsection{Condition for real eigenvalues}

As we saw, all the eigenvalues of $B = I - J^*$ need to be real if we 
apply the Chebyshev inertial iteration to a fixed-point iteration. 

The Jacobian of $f(x) = g(Ax+ b)$ at the fixed point $x^*$ is given by 
$
	J^* = Q^* A,	
$
where  $Q^* = \{q^*_{i,j}\} \in \mathbb{R}^{n \times n}$
is a diagonal matrix with the diagonal elements
\begin{equation}
	q^*_{i,i} = g'\left(\sum_{j=1}^n a_{i,j} x_j^* + b_i \right),\ i = 1,2,\ldots,n.
\end{equation}
It should be remarked that $Q^*A$ is not necessarily symmetric even when 
$A$ is symmetric. However, all the eigenvalues of $J^* = Q^* A$ are real
if the conditions described in the following theorem are satisfied.

\begin{theorem} \label{real_th}
	Assume that $g: \mathbb{R} \rightarrow \mathbb{R}$ is a differentiable function 
	satisfying $g'(x) \ge 0$ in the domain of $g$.
	If $A$ is a symmetric full-rank matrix, then 
	each eigenvalue of $J^* = Q^* A$ is real.
\end{theorem}
Theorem \ref{real_th} provides a sufficient condition for $J^*$ to have real eigenvalues.

\subsection{Proximal gradient descent}

A proximal iteration including a gradient descent process 
is widely used for implementing proximal gradient methods and projected gradient 
methods \cite{Parikh14} \cite{Combettes11}.
The proximal operator $\text{prox}_{h}(\cdot): \mathbb{R}^n \rightarrow \mathbb{R}^n $ 
for the function $h: \mathbb{R}^n \rightarrow \mathbb{R}^n$ is defined by
$
	\text{prox}_{h}(v) := \arg \min_{x} \left( h(x) +  (1/2)\|x - v\|^2 \right).
$
Assume that we want to solve the following minimization problem:
\begin{equation}  \label{objective}
	\text{minimize} \ u(x) + v(x),
\end{equation}
where $u: \mathbb{R}^n \rightarrow \mathbb{R}$ and $v: \mathbb{R}^n \rightarrow \mathbb{R}$.
The function $u$ is differentiable and $v$ is a strictly convex function which is not necessarily differentiable.
The iteration of the proximal gradient method is given by
$
	x^{(k+1)} = \text{prox}_{\lambda v}\left( x^{(k)} - \lambda \nabla u(x^{(k)}) \right),
$
where $\lambda$ is a step size. It is known that $x^{(k)}$ converges to the minimizer of (\ref{objective})
if $\lambda$ is included in the semi-closed interval $(0, 1/L]$ where $L$ is the Lipschitz constant of $\nabla u$.

In the following part of this subsection, we will discuss the proximal iteration 
regarding a regularized least square (LS) problem.
Suppose that a symmetric positive definite matrix $M \in \mathbb{R}^{n \times n}$ and
$c \in \mathbb{R}^n$ are given. 
The objective function of a regularized LS problem can be summarized as 
$
	U(x) + V(x) := x^T M x - c^T x + \lambda v(x),
$
where $U(x)$ is the quadratic function $x^T M x - c^T x$ and
$V(x)$ is a convex function (not necessarily differentiable) $\lambda v(x)$.
The $\lambda v(x)$ is a regularization term and the constant $\lambda$ is the regularization constant.

In the following, the proximal operator $\text{prox}_{\lambda v}(\cdot)$ is assumed to be 
a component-wise function denoted by $g(\cdot) := \text{prox}_{\lambda v}(\cdot)$.
Since the gradient of $U(x)$ is $\nabla U(x) = M x^{(k)} - c$,
the proximal gradient iteration for minimizing $U(x) + V(x)$ is given by 
\begin{equation}
	x^{(k+1)} = g\left((I - \gamma M) x^{(k)} + \gamma c \right),
\end{equation}
where $\gamma$ is the step size parameter.
Let $A :=  I - \gamma M$ and  $b := \gamma c$.
The proximal gradient descent process can be recast as a proximal iteration.
Note that the matrix $A$ becomes a symmetric matrix 
because $I - \gamma M$ is symmetric.

\section{Related Works}

The SOR method for linear equation solvers \cite{Young71} is useful especially for 
solving sparse linear equations. The convergence analysis can be found in \cite{Saad03}.
For acceleration of gradient descent (GD) methods, the heavy ball method involving 
an inertial or momentum term was proposed by Polyak \cite{Polyak64}.
Nesterov's accelerated GD for convex problems is another accelerated GD \cite{Nesterov83} with an inertial term.
The Chebyshev semi-iterative method \cite{Golub89, Mason03} is a fast method for solving linear equations 
that employs Chebyshev polynomials for bounding the spectral radius of an iterative system.
By using the recursive definition of the Chebyshev polynomials, one can derive
a momentum method based on the Chebyshev polynomial \cite{Golub89} and 
it provides the optimal linear convergence rate among the first order method.
It can be remarked that the Chebyshev semi-iterative method cannot be applied to 
nonlinear fixed-point iterations treated in this paper.

The Krasnosel'ski\v{\i}-Mann iteration \cite{Bauschke11b} is a method similar to the inertial iterations.
For a non-expansive operator $T$, the fixed-point iteration of the  Krasnosel'ski\v{\i}-Mann iteration is given as
$
	x^{(k+1)} = x^{(k)} + \lambda^{(k)} (T x^{(k)} - x^{(k)}).
$
It is proved that the point $x^{(k)}$ converges to a point in the set of fixed points of $T$
if $\lambda^{(k)} \in [0,1]$ satisfies $\sum_k \lambda^{(k)} (1 - \lambda^{(k)}) = +\infty$.

Acceleration of proximal gradient methods is a hot topic for convex optimization \cite{Bubeck15}.
FISTA \cite{Beck09} is an accelerated algorithm based on ISTA and Nesterov method \cite{Nesterov83}.
FISTA and its variants \cite{Bioucas-Dias07}
achieves much faster global convergence speed $O(1/k^2)$ compared with that of the original ISTA, i.e., $O(1/k)$.

\section{Experiments}

\subsection{Linear fixed-point iteration: Jacobi method}

\begin{figure}[htbp]
\begin{center}
\centerline{\includegraphics[width=\columnwidth]{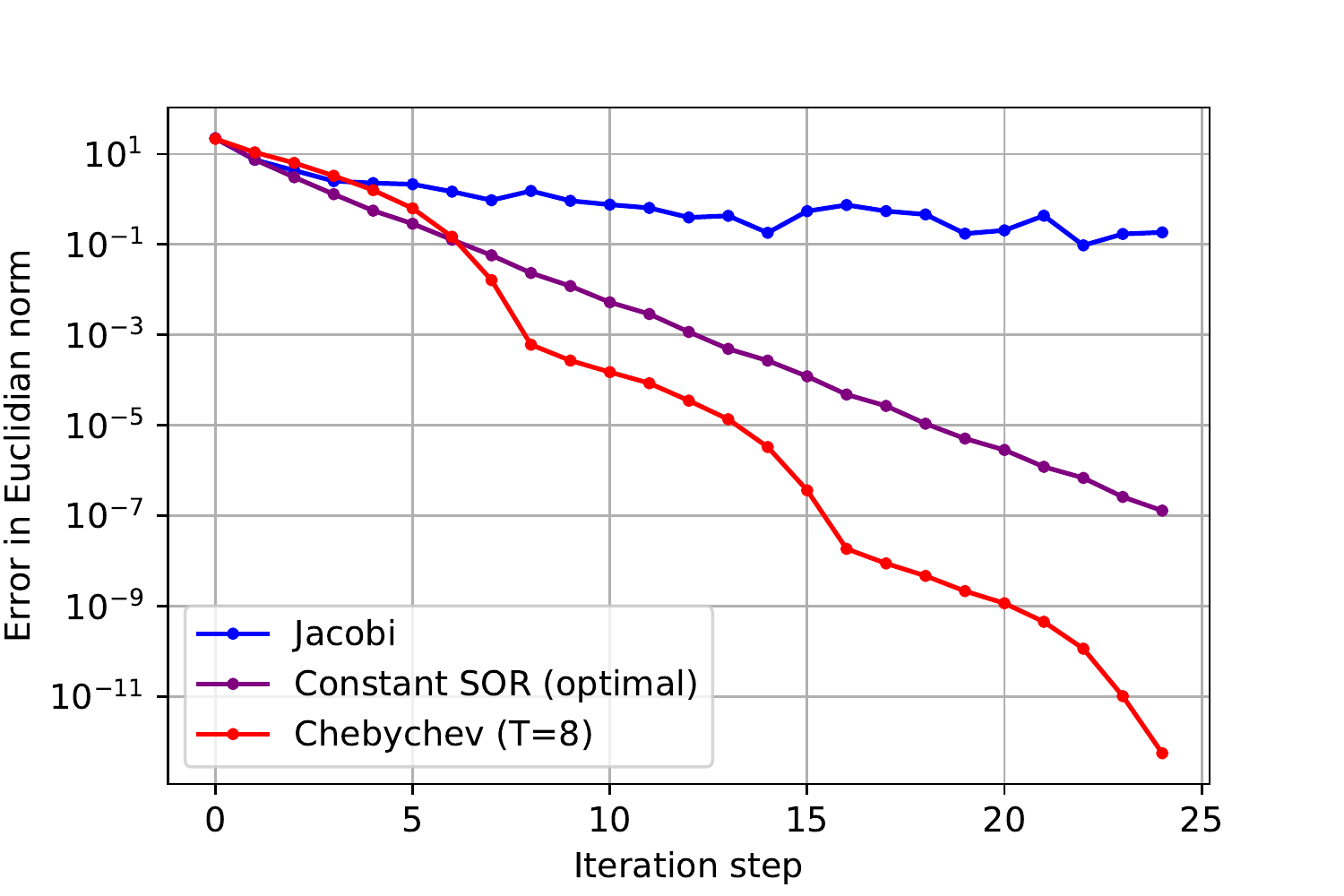}} 
\caption{Jacobi method, Constant SOR factor, and Chebyshev inertial iteration: 
The linear equation $P x = 0$ is assumed without loss of generality where
$P = I + M^T M \in \mathbb{R}^{512 \times 512}$ is a diagonally dominant matrix.
Each element of $M$ follows ${\cal N}(0, 0.03^2)$. The minimum and maximum 
eigenvalues of $B = I-A$ are given by 
$\lambda_{max}(B) = 1.922$ and $\lambda_{min}(B) = 0.6766$, respectively.
The 3 curves represent the errors in Euclidean norm from the fixed point, i.e., $\|x^{(k)} - x^*\|$ where
$x^*  = 0$.
The optimal constant SOR factor 
is given by $\omega_{SOR} = 2/(\lambda_{min}(B) + \lambda_{max}(B)) = 0.7697$ \cite{Saad03}.
}
\label{jacobi}
\end{center}
\vskip -0.2in
\end{figure}

The Jacobi method is a linear fixed-point iteration
\begin{equation}
	x^{(k+1)} = g(A x^{(k)} + b)	 = - D^{-1} (P - D) x^{(k)} + D^{-1} q 
\end{equation}
for solving linear equation $P x = q$ where
$A := - D^{-1} (P - D)$, $b := D^{-1} q$, and $g(x) := x$.
We can apply the Chebyshev inertial iteration to the Jacobi method for accelerating the convergence.
Figure \ref{jacobi} shows the convergence behaviors.
It should be remarked that the optimal constant SOR factor \cite{Saad03} is 
exactly the same as the Chebyshev inertial factor with $T=1$.
We can see that the Chebyshev inertial iteration provides much faster convergence 
compared with the plain Jacobi and the constant SOR factor method.
The error curve of the Chebyshev inertial iteration indicates a wave-like shape.
This is because the error is tightly bounded periodically as shown in Theorem \ref{convergence_rate}, i.e., 
the error is guaranteed to be small when the iteration index is a multiple of $T$.

\subsection{Nonlinear fixed-point iterations}

\begin{figure}[htbp]
\vskip 0.2in
\begin{center}
\centerline{\includegraphics[width=\columnwidth]{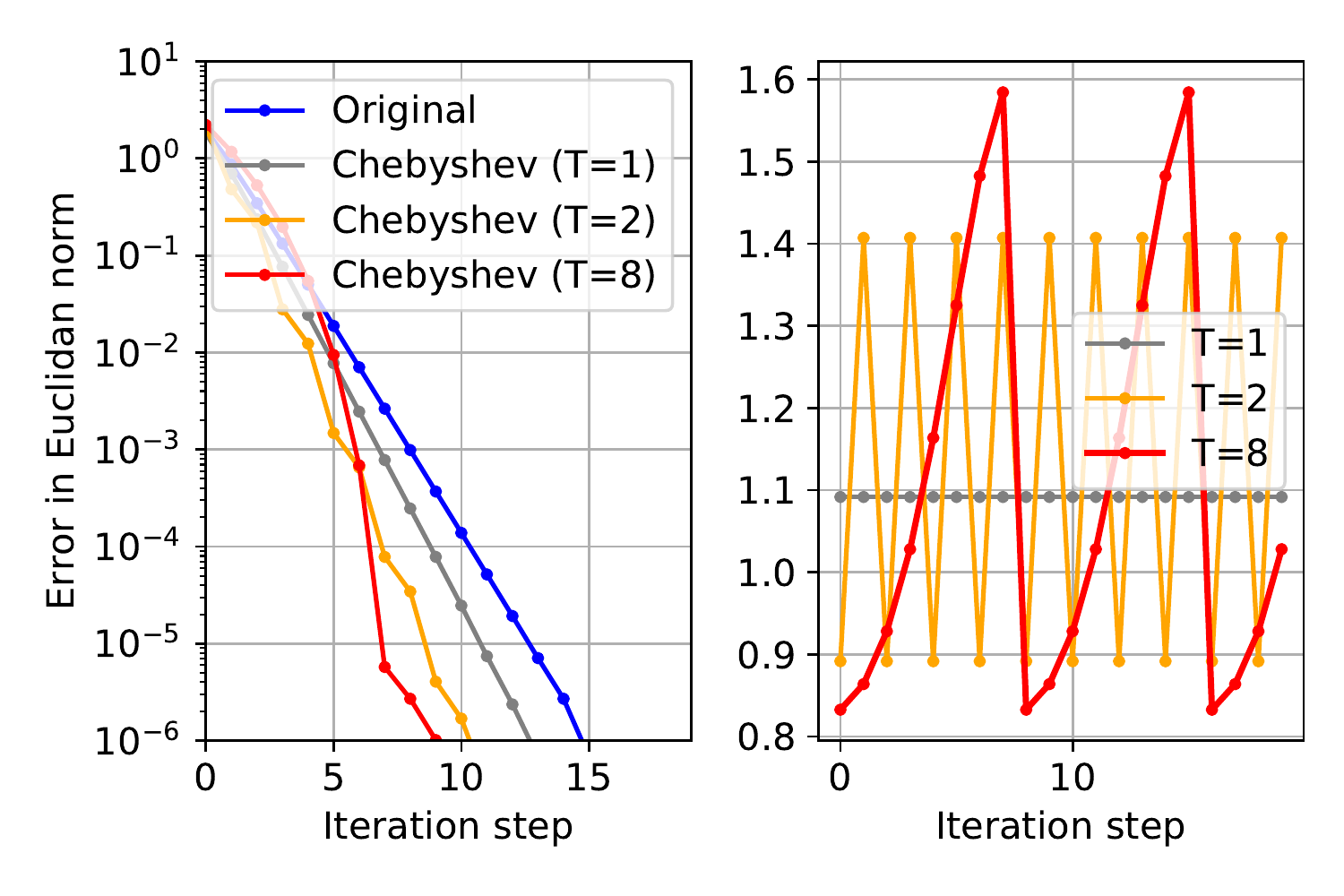}}
\caption{Nonlinear fixed-point iteration: 
The left figure shows the error $\|x^{(k)} - x^*\|$ for 
the original fixed point iteration and the Chebyshev inertial iterations $(T= 1,2,8)$.
The fixed point  is $x^* = (2.96, 2.96)^T$.
The minimum and maximum eigenvalues of $B = I- J^*$  are 
$\lambda_{max}(B) = 1.216$, $\lambda_{min}(B) = 0.626$, respectively.
The horizontal axes represents iteration step. The right figure shows
the Chebyshev inertial factors for $T = 1, 2, 8$ as a function of the iteration step.
}
\label{polynomial}
\end{center}
\vskip -0.2in
\end{figure}

A nonlinear function $f = (f_1, f_2)^T: \mathbb{R}^2 \rightarrow \mathbb{R}^2$
\begin{eqnarray}
f_1(x_1, x_2) &:=& x_1^{0.2} + x_2^{0.5}, \\
f_2(x_1, x_2) &:=& x_1^{0.5} + x_2^{0.2}
\end{eqnarray}
is assumed here.
Figure \ref{polynomial} (left) shows the error curves as functions of
iteration. We can observe that the Chebyshev inertial iteration actually accelerates the convergence to the
fixed point. The Chebyshev inertial iteration results in a zigzag shape as 
depicted in Fig. \ref{polynomial} (right).

\begin{figure}[htbp]
\vskip 0.2in
\begin{center}
\centerline{\includegraphics[width=\columnwidth]{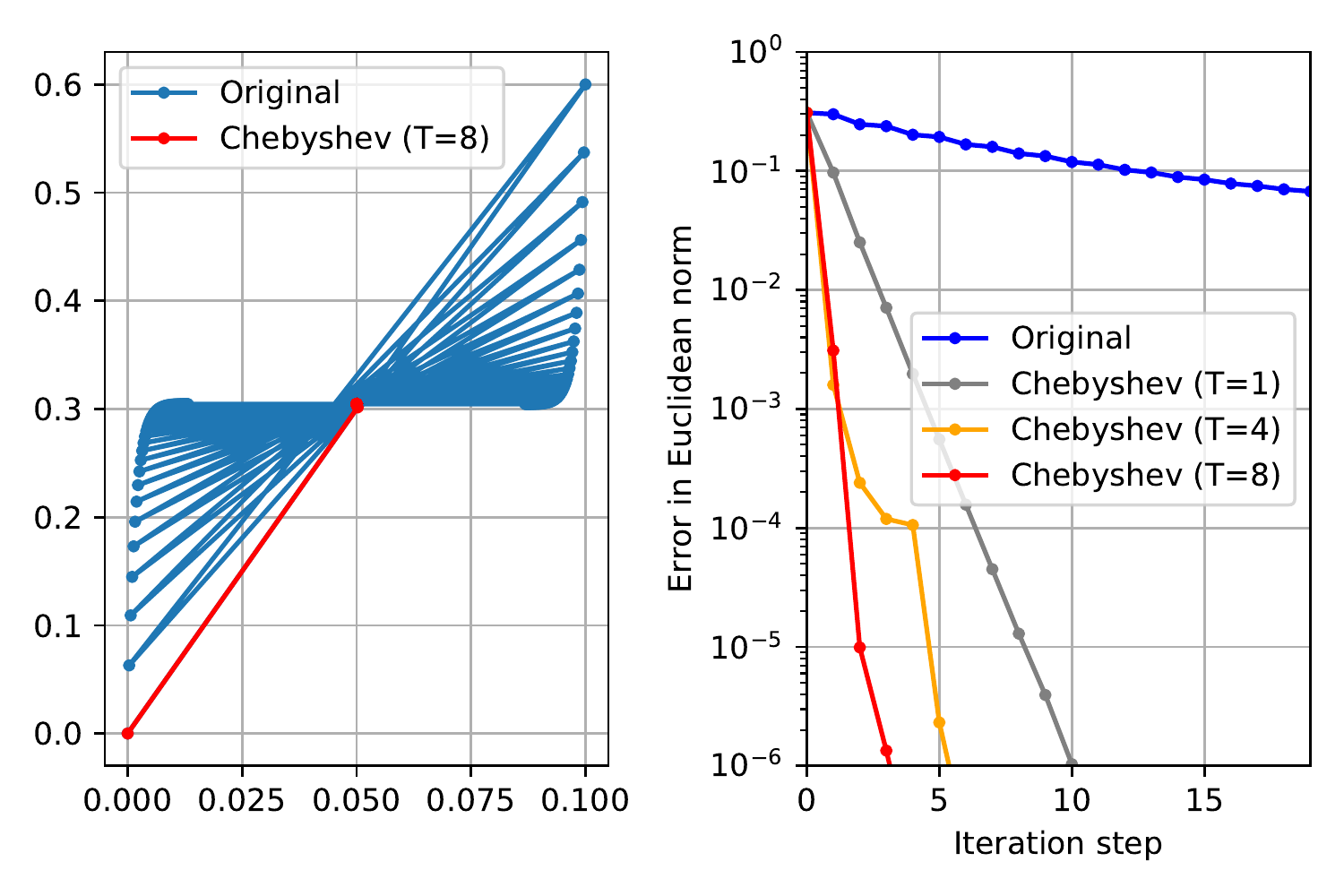}}
\caption{Solving a nonlinear equation: The solution of a nonlinear equation $y = x + \tanh (x) (y = (0.1, 0.6)^T)$ 
can be found by the fixed point iteration $x^{(k+1)} = y - \tanh(x^{(k)})$. 
The fixed point is $(0.0500, 0.3045)^T$.
The left figure indicates the search trajectories of the original fixed-point iteration and 
the accelerated iteration by the Chebyshev inertial iteration. 
The right figure presents the  error from the fixed point. }
\label{nonlinear}
\end{center}
\vskip -0.2in
\end{figure}

If the initial point is sufficiently close to the solution, 
one can solve  a nonlinear equation with a fixed-point iteration.
Figure \ref{nonlinear} explains the behaviors of two-dimensional fixed-point iterations for solving
a non-linear equation $y = x + \tanh(x)$. We can see that the Chebyshev inertial  provides 
much steeper error curves compared with the original iteration in Fig. \ref{nonlinear} (right).

\begin{figure}[htbp]
\vskip 0.2in
\begin{center}
\centerline{\includegraphics[width=\columnwidth]{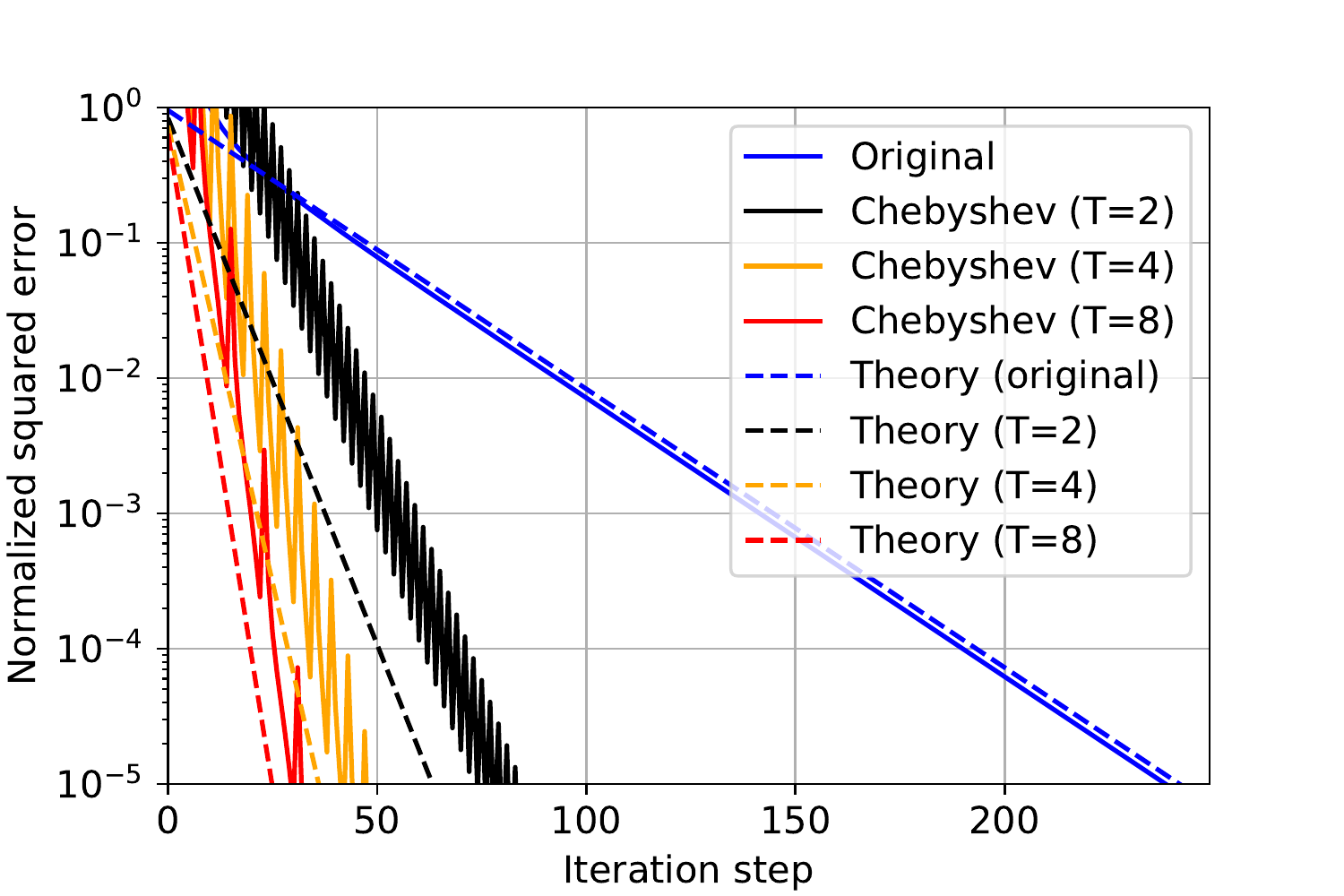}}
\caption{Nonlinear fixed-point iteration in a high dimensional space: 
The fixed-point iteration is $x^{(k+1)} = \tanh(A x^{(k)})$ and $A$ is a random Gram matrix 
$A = M^T M$ of size $512 \times 512$.
Each element of $M$ follows i.i.d. Gaussian distribution with mean zero and 
standard deviation 0.022. 
 The solid lines represent the normalized squared error $\|x^{(k)} - x^*\|^2/n$ for 
the original fixed point iteration and the Chebyshev inertial iterations $(T= 2,4,8)$.
The fixed point is the origin. The doted lines (Theory) indicate $\rho(J^*)^{2k}$ 
(original),
and $[q_{CI}(T)]^{2k}$ 
where $T = 2, 4, 8$, and $a = 0.0234, b = 1.0$.
}
\label{theory}
\end{center}
\vskip -0.2in
\end{figure}

Consider the proximal iteration
$
	x^{(k+1)} = \tanh(A x^{(k)})		
$
where  $A \in \mathbb{R}^{512 \times 512}$.
The Jacobian $J^*$ at the fixed point is identical 
to $A$ and the minimum and maximum eigenvalues of $J^*$ are $0$ and $0.9766$, respectively.
Figure \ref{theory} shows the normalized errors $\|x^{(k)} - x^*\|/n$ for the 
original and the Chebyshev inertial iterations. This figure also includes 
the theoretical values $[\rho(J^*)]^{2k}$ and 
$[q_{CI}(T)]^{2 k}$ 
under the condition $a = 0.0234$ and $b = 1.0$.
We can observe that the error curves of the Chebyshev inertial iteration shows a zigzag shape as well.
The lower envelopes of the zigzag lines (corresponding to periodically bounded errors)
and  the theoretical values derived from Theorem \ref{convergence_rate}
have almost identical slopes. 

\subsection{ISTA for sparse signal recovery}

ISTA \cite{Daubechies04}
 is a proximal gradient descent method designed for sparse signal recovery problems.
The problem setup of the sparse signal recovery problem discussed here is as follows.
Let $x \in \mathbb{R}^n$ be a source sparse signal. We assume each element in $x  = (x_1,\ldots, x_n)^T$ follows
i.i.d. Berrnoulli-Gaussian distribution, i.e, non-zero elements in $x$ occur with probability $p$
and these non-zero elements follow the normal distribution ${\cal N}(0, 1)$.
Each element in a sensing matrix $M \in \mathbb{R}^{m \times n}$ is assumed to be generated
randomly according to  ${\cal N}(0, 1)$. The observation signal $y$ is given by
$
	y = M x + w,
$
where $w \in \mathbb{R}^m$ is an i.i.d. noise vector whose elements follow  ${\cal N}(0, \sigma^2)$.
Our task is to recover $x$ from $y$ as correct as possible.
This problem setting is also known as compressed sensing \cite{Donoho06, Candes06}.

A common approach to tackle the sparse signal recovery problem described above is 
to use Lasso formulation \cite{Tibshirani96, Efron04}:
\begin{equation}
\hat x := \arg \min_{x \in \mathbb{R}^n} \frac 1 2 \| y - M x\|^2 + \lambda \|x\|_1.	
\end{equation}
ISTA is the proximal gradient descent method for minimizing the Lasso objective function.
The fixed-point iteration of ISTA is given by
\begin{eqnarray} \label{ISTAiteration1}
	r_t &=& s_t + \beta M^T ( y - M s_t), \\ \label{ISTAiteration2}
	s_{t+1} &=& \eta(r_t; \tau),
\end{eqnarray}
where $\eta(x; \tau)$ is the soft shrinkage function defined by
$
	\eta(x; \tau) := \text{sign}(x) \max \{|x| - \tau, 0 \}.	
$

\begin{figure}[htbp]
\begin{center}
\centerline{\includegraphics[width=\columnwidth]{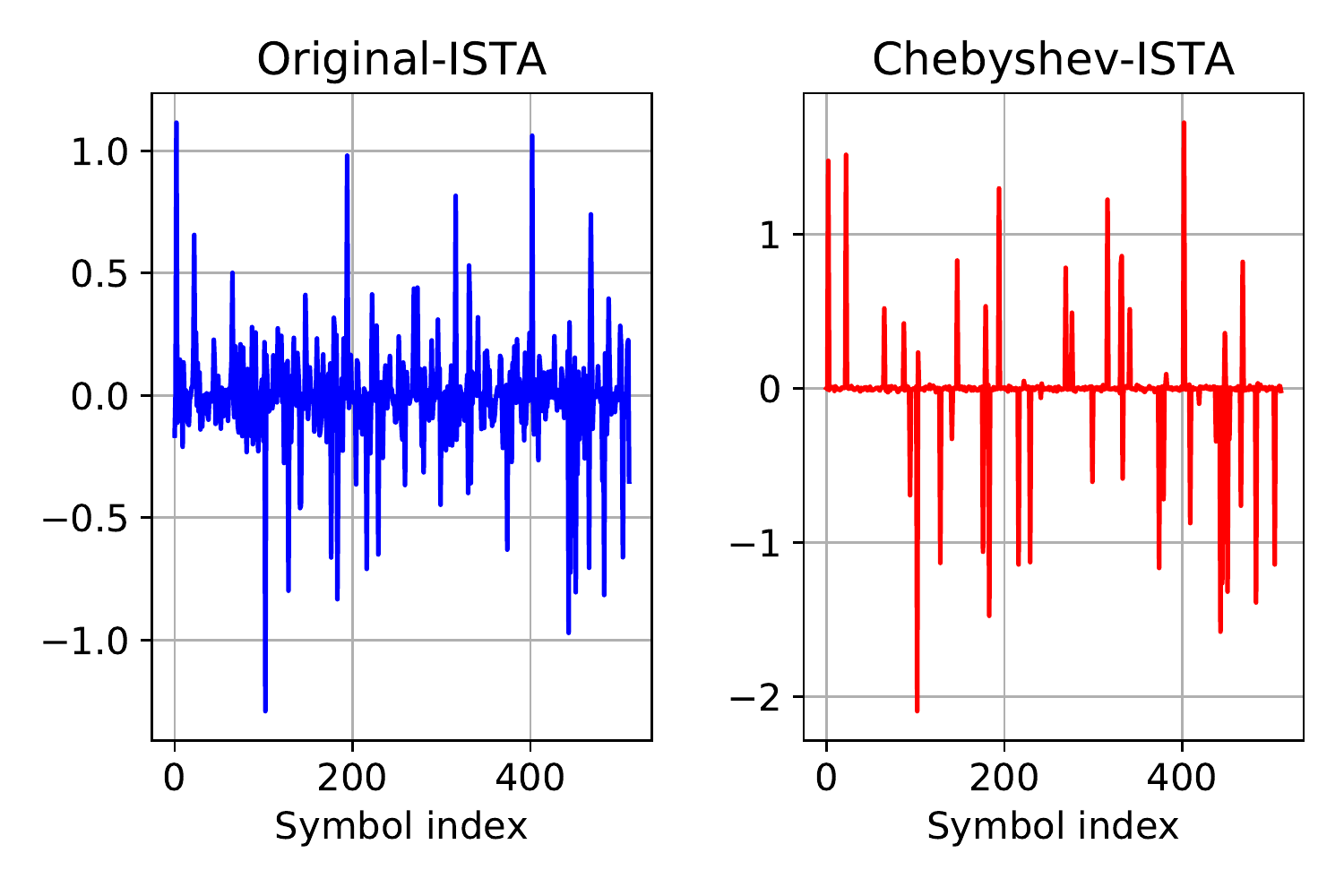}}
\caption{Reconstructed signals by ISTA: The system parameters are $n = 512, m = 256, p = 0.1, \sigma = 0.1$. 
	The left and right figures indicate the recovered signal after 200 iterations of
	ISTA and Chebyshev-ISTA, respectively. The source sparse signal $x$ is the same 
	in both figures.
}
\label{ista_12}
\end{center}
\vskip -0.2in
\end{figure}

The softplus function $s_p(x)$ is defined as 
\begin{equation}
	s_p(x) := ({1}/{\beta})\log(1  + \exp(\beta x)).		
\end{equation}
In the following experiments, we will use 
a {\em differentiable soft shrinkage function} defined by
\begin{equation}\label{eta}
	\tilde \eta(x; \tau) := s_p(x - \tau) + s_p(- (x + \tau))	
\end{equation}
with $\beta = 100$ instead of the soft shrinkage function $\eta(x; \tau)$ because
$\tilde \eta(x; \tau) $ is differentiable.

Let $G = M^T M$. The step size parameter and shrinkage parameter are 
set to $\gamma = \tau = 1/\lambda_{max}(G)$.
We are now ready to write 
the ISTA iteration (\ref{ISTAiteration1}), (\ref{ISTAiteration2}) in the form of proximal iteration:
\begin{equation}
	x^{(k+1)} = \tilde \eta(A x^{(k)} + b ; \tau),
\end{equation}
where 
$A := I - \gamma G$ and $b := \gamma M^T y$.

Figure \ref{ista_12} presents recovered signals by the original ISTA and Chebyshev-ISTA.
The reconstructed signals are the results obtained after
200 iterations. It can be immediately observed that Chebyshev-ISTA provides
sparser reconstructed signal than that of the original ISTA. This is because 
the number of iterations (200 iterations) is not enough for the original ISTA 
to achieve reasonable reconstruction.

\begin{figure}[htbp]
\begin{center}
\centerline{\includegraphics[width=\columnwidth]{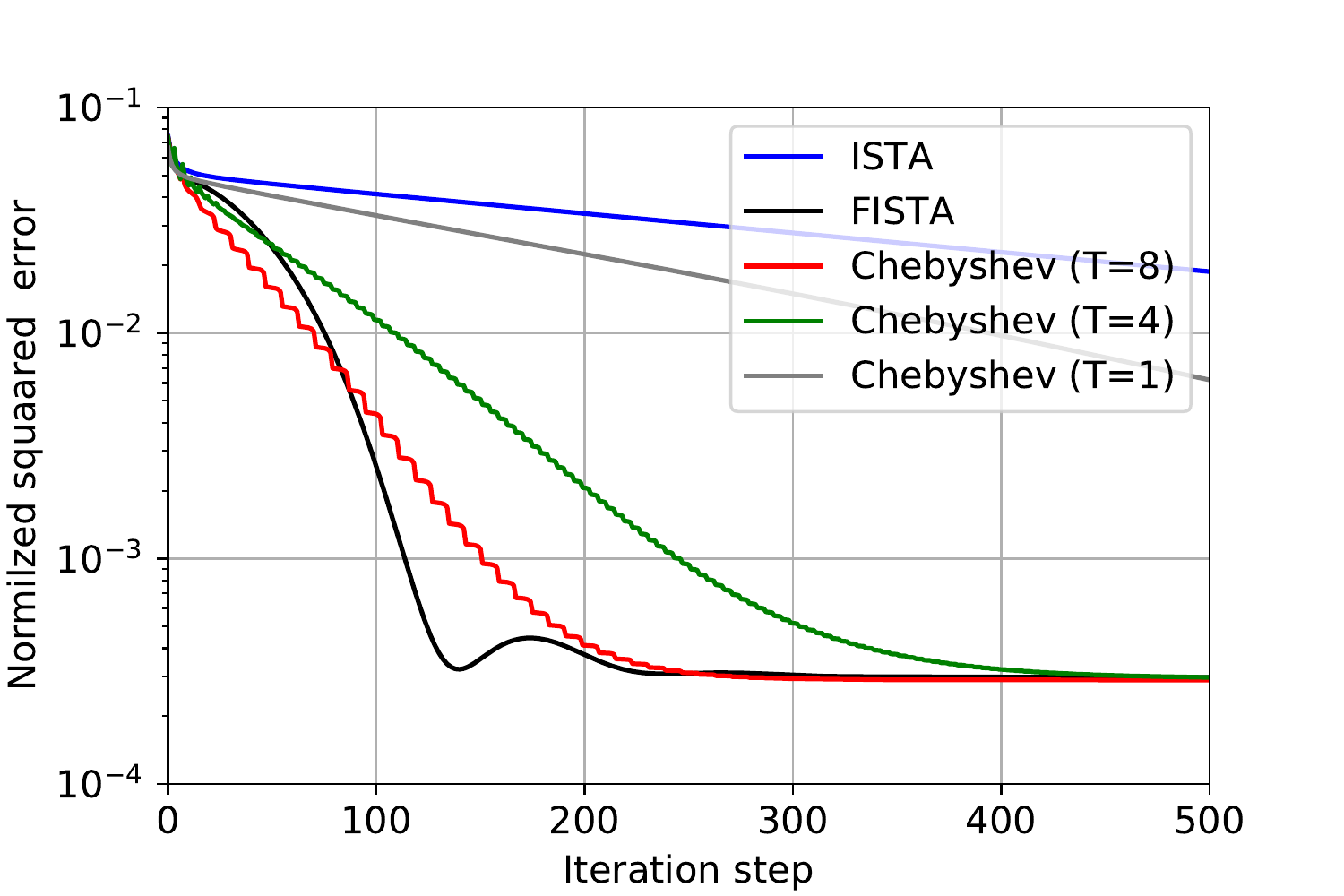}}
\caption{Sparse signal recovery by ISTA: $n = 512, m = 256, p = 0.1, \sigma = 0.1$. 
The observation vector is modeled by $y = M x + w$ where $M \in \mathbb{R}^{m \times n}$.
The normalized squared error $\|x^{(k)} - x\|^2/n$ is averaged over 1000 trials. As a reference, 
the error curve  of FISTA \cite{Beck09} is also included.}
\label{ista_average}
\end{center}
\vskip -0.2in
\end{figure}

Figure \ref{ista_average} indicates the averaged error for 1000 trials.
The Chebyshev-inertial ISTA ($T=8$) requires only 250--300 iterations 
to reach the point where the original ISTA achieves with 3000 iterations.
The convergence speed of the Chebyshev-inertial 
ISTA is almost comparable to that of FISTA but 
the proposed scheme provides smaller squared errors 
when the number of iterations is less than 70.
This result implies that the proposed scheme may be used 
for convex or non-convex problems 
as an alternative to FISTA and the Nesterov's method.

\subsection{Modified Richardson iteration for deblurring}

\begin{figure}[htbp]
\begin{center}
\centerline{\includegraphics[width=\columnwidth]{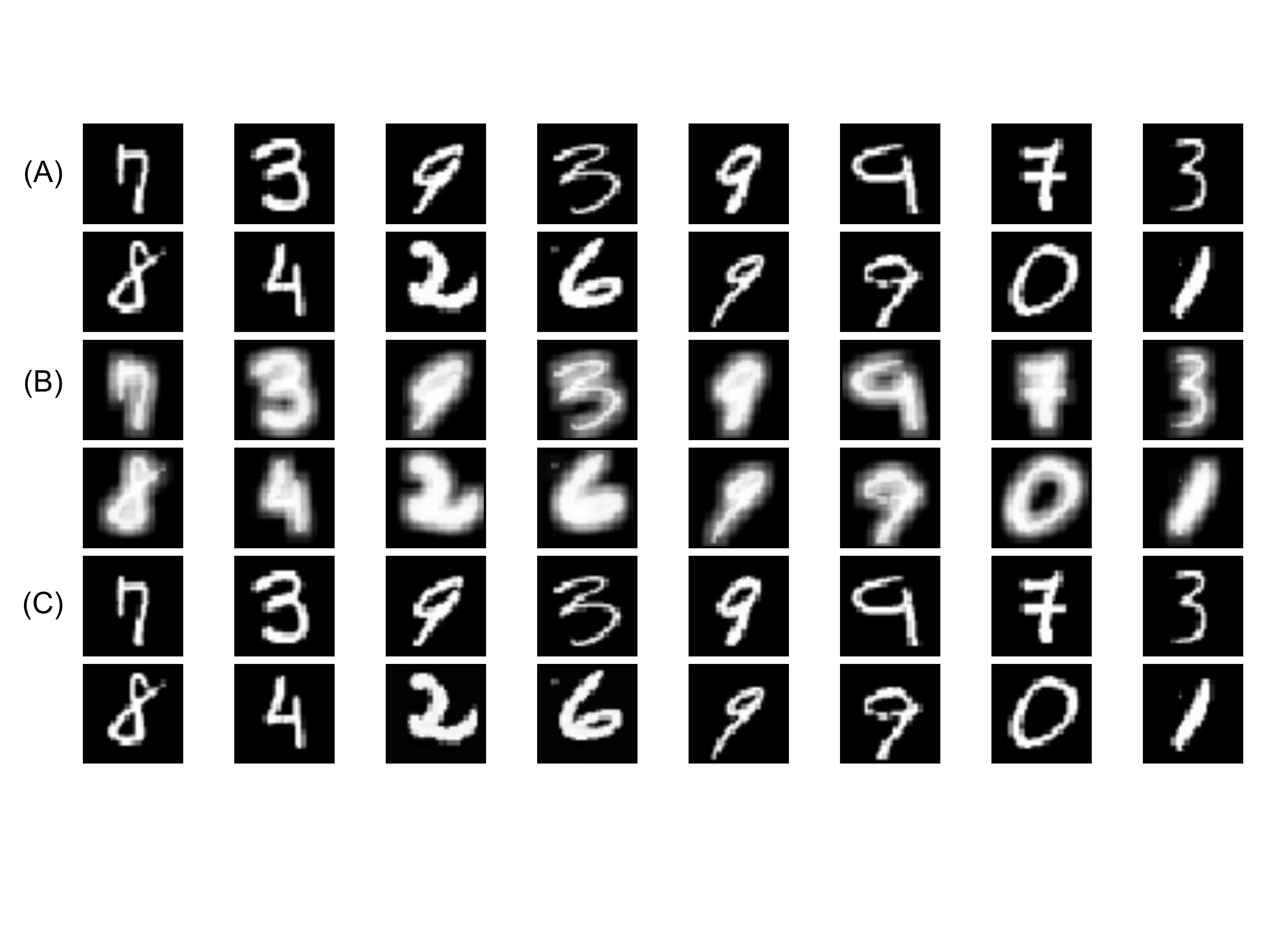}}
\vspace{-1.0cm}
\centerline{\includegraphics[width=\columnwidth]{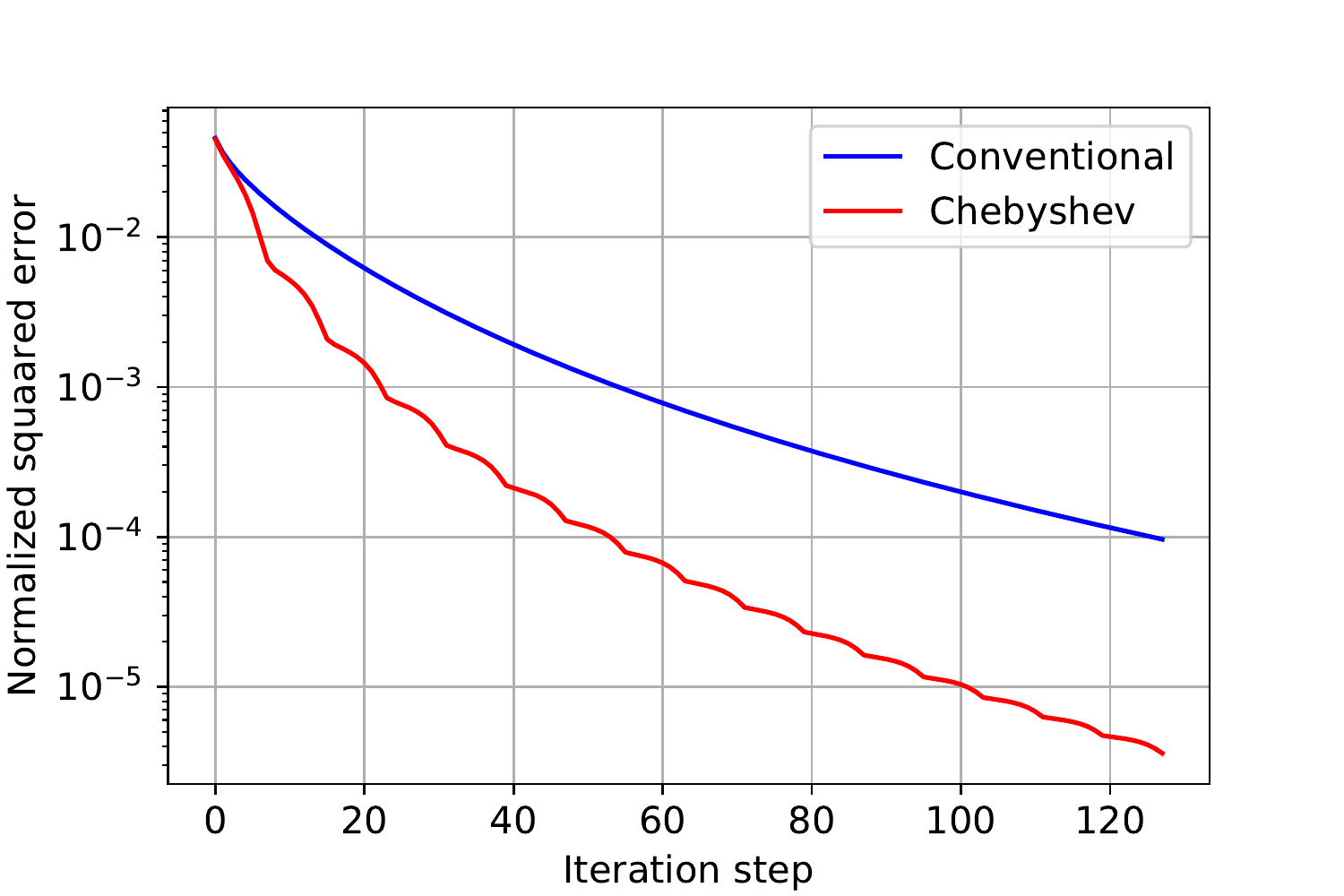}}
\caption{Nonlinear deblurring by modified Richardson iteration: The blurring process is the following: A $7 \times 7$ convolutional kernel with the value 0.1 except for the center (the value 1.5) 
 is first applied to a normalized MNIST image in $[0,1]^{28 \times 28}$ and then a sigmoid function is applied to the convoluted image. The parameter $\omega$ is set to $0.8$.
(A) Original image, (B) Blurred image, (C) Recovered image by Richardson iteration with the Chebyshev inertial iteration  (number of iterations is 128, $T=8$, $\lambda_{min}(B) = 0.18, \lambda_{max}(B) = 0.98$).
The graph represents the normalized errors $\|x^{(k)} - x\|^2/28^2$ of the original fixed-point iteration and the accelerated iteration.}
\label{blur}
\end{center}
\vskip -0.2in
\end{figure}

Assume that a signal $x \in \mathbb{R}^n$ is distorted by a nonlinear process $y = f(x)$.
The modified Richardson iteration \cite{Saad03} is usually known as a method for
solving a linear equation but it can be applied to solve a nonlinear equation.
The fixed-point iteration is given by 
\begin{equation}
	x^{(k+1)} = x^{(k)} + \omega ( y - f(x^{(k)} )  ),
\end{equation}
where $\omega$ is a positive constant. 
Figure \ref{blur} shows the case where $f$ expresses a nonlinear blurring process to the MNIST images.
We can see that the Chebyshev inertial iteration can accelerate the convergence as well.

\section{Conclusion}

A novel method for accelerating the convergence speed of 
fixed-point iterations is presented in this paper.
The analysis based on linearization around a fixed point 
unveils the local convergence behavior of the Chebyshev inertial iteration.
Experimental results shown in the paper support the theoretical analysis on the Chebyshev inertial iteration.
For example, the results shown in Fig. \ref{theory} indicate excellent agreement 
with the experimental results.
The proposed method appears effective in particular for accelerating the proximal gradient methods 
such as ISTA.  The proposed scheme requires no training process to adjust the inertial factors 
compared with the deep-unfolding approach \cite{Gregor10, Ito19} for accelerating ISTA.
The Chebyshev inertial iteration can be applied to many iterative algorithms with fixed-point iterations 
without requiring additional computational cost. 

\section*{Acknowledgement}
This work was supported by JSPS Grant-in-Aid for Scientific Research (B) 
Grant Number 16H02878.

\bibliography{fixedpoint}
\bibliographystyle{ieeepes}

\appendix

\subsection{Proof of Lemma 1}

Since any eigenvalue of $J^*$ is real, 
all the eigenvalues of $B = I - J^*$ are also real.
Because $f$ is a locally contracting mapping,
$\rho(J^*) < 1$ should be satisfied. 
 This means that
	$\lambda_{min}(J^*) >  -1$ and $\lambda_{max}(J^*) <  1$. It is thus clear that
	$\lambda_{min}(B) >  0$ and $\lambda_{max}(B) <  2$ are satisfied.

\subsection{Proof of Lemma 2}
From the definition of the ATMC polynomial, we have 
\begin{eqnarray}
\tilde \beta_T(\lambda) 
&=&  \left[\hat C_T(0; a, b) \right]^{-1}
\hat C_T(\lambda; a, b) \\
&=& \frac{C_T(\lambda; a, b)}{C_T(0; a, b)} \\
&=& \frac{C_T\left(\frac{2 \lambda - b - a}{b- a} \right)}{C_T\left(- \frac{b + a}{b- a} \right)}.
\end{eqnarray}
Due to the property of the Chebyshev polynomial such that $|C_T(x)| \le 1$ for $x \in [-1,1]$,
the absolute value of the numerator is bounded as
\begin{equation}
\left|C_T\left(\frac{2 \lambda - b - a}{b- a} \right) \right|  \le 1
\end{equation}
for $a \le \lambda \le b$. By using this inequality,  the absolute value of $\tilde \beta_T(\lambda)$ 
can be upper bounded by
\begin{eqnarray}
|\tilde \beta_T(\lambda) |
&\le&   \left|C_T\left(- \frac{b+a}{b-a} \right)\right| 
^{-1}\end{eqnarray}
for $a \le \lambda \le b$.

It is known that the Chebyshev polynomial can be expressed as 
$
C_k(x) = (-1	)^k \cosh(k \cosh^{-1}(-x))
$
if $x \le -1$ \cite{Mason03}. By using this identity, 
the absolute value of $\tilde \beta_T(\lambda)$ can be upper bounded by
\begin{eqnarray}  \nonumber
|\tilde \beta_T(\lambda)| &\le& \left|C_T\left(- \frac{b+a}{b-a} \right)\right| 
^{-1} \\
&=& \text{\rm sech}\left(T \cosh^{-1}\left(\frac{b+a}{b-a}\right) \right),
\end{eqnarray}
where $1/\cosh(x) = \text{\rm sech}(x)$.

\subsection{Proof of Theorem 1}
Expanding $S^{(k)}(x)$ around the fixed point $x^*$, we have 
\begin{equation}
	S^{(k)}(x) = S^{(k)}(x^*) +  \left(I - \omega_{k} B \right) (x - x^*) 
	+ o(\|x - x^*\|).
\end{equation}
By substituting $x = x^{(k)}$,  the above equation can be rewritten as
\begin{equation}
	x^{(k+1)} - x^* =  \left(I - \omega_{k} B \right) (x^{(k)} - x^*) + o(\|x^{(k)} - x^*\|).
\end{equation}

Let
\begin{equation}
\xi := \max_{0 \le k \le T-1} \frac{\|x^{(k)} - x^*\|}{\|x^{(0)} - x^*\|}
\end{equation}
and $\delta := \xi \|x^{(0)} - x^*\|$.  Note that $\xi$ is finitely bounded 
because of the definition of the Chebyshev inertial factors.
From these definitions,
we have
\begin{equation}
	\frac{o(\|x^{(k)} - x^*\|) }{\delta} \rightarrow 0.
\end{equation}
for $k = 0, 1, \ldots, T-1$ if $\delta \rightarrow 0$. 
Replacing $o(\|x^{(k)} - x^*\|)$ by $o(\delta)$, we obtain
\begin{equation}
	x^{(k)} - x^* = \left(I - \omega_{k-1} B \right) (x^{(k-1)} - x^*) + o(\delta)
\end{equation}
for $k = 1, 2, \ldots, T$.
This equation leads to
\begin{equation}
	x^{(T)} - x^* =
\left(\prod_{k = 0}^{T-1} \left(I - \omega_{k} B \right)  \right) (x^{(0)} - x^*) + o(\delta).
\end{equation}
By taking the norm of the both sides, we have
\begin{eqnarray} \nonumber
	\|x^{(T)} - x^*\| &=&
\left\|\left(\prod_{k = 0}^{T-1} \left(I - \omega_{k} B \right)  \right) (x^{(0)} - x^*) + o(\delta) \right\| \\  \nonumber
&\le&
\left\|\left(\prod_{k = 0}^{T-1} \left(I - \omega_{k} B \right)  \right) (x^{(0)} - x^*)\right\| + \left\|o(\delta) \right\| \\ \nonumber 
&\le& 
\left\|\left(\prod_{k = 0}^{T-1} \left(I - \omega_{k} B \right)  \right) \right\| 
\left\|(x^{(0)} - x^*)\right\| + \left\|o(\delta) \right\| \\ \nonumber
&=& 
\rho \left(\prod_{k = 0}^{T-1} \left(I - \omega_{k} B \right)  \right)
\left\|(x^{(0)} - x^*)\right\| + \left\|o(\delta) \right\|, \\
\end{eqnarray}
where the first inequality is due to the triangle inequality and
the second inequality is due to sub-multiplicativity of the operator norm.
By dividing both sides by $\|x^{(0)} - x^*\|$, we get
\begin{equation}
	 \lim_{\delta \rightarrow 0}\frac{\|x^{(T)} - x^*\|}{\|x^{(0)} - x^*\|}
	 \le \rho \left(\prod_{k = 0}^{T-1} \left(I - \omega_{k} B \right)  \right).
\end{equation}
Since the inertial factors $\omega_k$ are periodical, we can apply the same argument for $k > T$
and obtain 
\begin{equation} \label{adjacent_ineq}
	 \lim_{\delta \rightarrow 0}\frac{\|x^{(\ell T)} - x^*\|}{\|x^{((\ell - 1)T)} - x^*\|}
	 \le \rho \left(\prod_{k = 0}^{T-1} \left(I - \omega_{k} B \right)  \right),
\end{equation}
where $\ell$ is a positive integer.
By using the inequality (\ref{adjacent_ineq}),  we have
\begin{eqnarray} \nonumber
&& \hspace{-0.8cm} \lim_{\delta \rightarrow 0}\frac{\|x^{(\ell T)} - x^*\|}{\|x^{(0)} - x^*\|}\\ \nonumber
&=& 	
\lim_{\delta \rightarrow 0}\left(
\frac{\|x^{(\ell T)} - x^*\|}{\|x^{((\ell-1) T)} - x^*\|}
\frac{\|x^{((\ell-1) T)} - x^*\|}{\|x^{((\ell-2) T)} - x^*\|}
\cdots \right)
\\
&\le& 
\rho\left(\prod_{k = 0}^{T-1} \left(I - \omega_k B \right)  \right)^{\ell}.
\end{eqnarray}
From this inequality, we immediately obtain the
claim of the theorem.

\subsection{Proof of Theorem 2}

From the assumption $g'(x) \ge 0$, $Q^*$ is a diagonal matrix with non-negative diagonal elements.
Since $A \in \mathbb{R}^{n \times n}$ is a full-rank matrix, the rank of $Q^* A$ coincides with
the number of non-zero diagonal elements in $Q^*$ which is denoted by $r$.

Darazin-Hynsworth theorem \cite{Drazin62} presents that the necessary and sufficient condition 
for $X \in \mathbb{C}^{n \times n}$ having $m (\le n)$ linearly independent eigenvectors 
and corresponding $m$ real eigenvalues is that there exists a symmetric semi-positive definite matrix $S$ 
with rank $m$ satisfying 
$
	X S = S X^H.
$
The matrix $X^H$ represents Hermitian transpose of $X$.
As we saw, $Q^*$ is a symmetric semi-positive definite matrix with rank $r$.
Multiplying $Q^*$ to $Q^* A$ from the right, we have the equality
\begin{equation}
	(Q^* A) Q^* = Q^* (A Q^*) =  Q^* (Q^{*T} A^T)^T = Q^* (Q^* A)^T,
\end{equation}
which satisfies the necessary and sufficient condition of Darazin-Hynsworth Theorem.
This implies that $Q^* A$ has $r$ real eigenvalues.

\end{document}